\documentclass[10pt]{amsart}
\usepackage{mathrsfs}
\usepackage{amsfonts} %%% i.e. use 12pt type
\usepackage{graphicx,latexsym,bm,amsmath,amssymb,verbatim,multicol,lscape}
\theoremstyle{definition}

\theoremstyle{remark}

\numberwithin{equation}{section}
\begin{document}
\title[On the periodicity of a class of arithmetic functions]
{On the periodicity of a class of arithmetic functions associated with multiplicative functions}%
\author{Guoyou Qian}
\address{Center for Combinatorics, Nankai University, Tianjin 300071, P.R. China}
\email{qiangy1230@gmail.com, qiangy1230@163.com}
\author{Qianrong Tan}
\address{School of Mathematics and Computer Science, Panzhihua University,
Panzhihua 617000, P.R. China} \email{tqrmei6@126.com}
\author{Shaofang Hong*}
\address{Yangtze Center of Mathematics, Sichuan University, Chengdu 610064, P.R. China and
 Mathematical College, Sichuan University, Chengdu 610064, P.R. China}
%\curraddr{}
\email{sfhong@scu.edu.cn, s-f.hong@tom.com, hongsf02@yahoo.com}
\thanks{*Hong is the corresponding author and was supported partially by the
National Science Foundation of China Grant \# 10971145 and by the Ph.D. Programs
Foundation of Ministry of Education of China Grant \#20100181110073}

%\subjclass{Primary 11T22,11R18}%
%\date{\today}%
%\dedicatory{}%
%\commby{}%
% ----------------------------------------------------------------
\begin{abstract}
Let $k\ge 1,a\ge 1,b\ge 0$ and $ c\ge 1$ be integers. Let $f$ be a
multiplicative function with $f(n)\ne 0$ for all positive integers $n$.
We define the arithmetic function $g_{k,f}$ for any positive integer
$n$ by $g_{k,f}(n):=\frac{\prod_{i=0}^k f(b+a(n+ic))} {f({\rm
lcm}_{0\le i\le k} \{b+a(n+ic)\})}$. We first show that $g_{k,f}$ is
periodic and $c {\rm lcm}(1,...,k)$ is its period. Consequently,
we provide a detailed local analysis to the periodic function
$g_{k,\varphi}$, and determine the smallest period of
$g_{k,\varphi}$, where $\varphi$ is the Euler phi function.
\end{abstract}

\maketitle

\section{\bf Introduction}
Chebyshev \cite{[Ch]} initiated the study of the least common
multiple of consecutive positive integers for the first significant
attempt to prove prime number theorem. An equivalent of prime
number theorem says that $\log {\rm lcm}(1, ...,n)\sim n$
as $n$ goes to infinity. Hanson \cite{[Ha]} and Nair \cite{[N]} got the
upper bound and lower bound of ${\rm lcm}_{1\le i\le n}\{i\}$
respectively. Bateman, Kalb and Stenger \cite{[BKS]} obtained an
asymptotic estimate for the least common multiple of arithmetic
progressions. Hong, Qian and Tan \cite{[HQT1]} obtained an
asymptotic estimate for the least common multiple of a sequence of
products of linear polynomials.

On the other hand, the study of periodic arithmetic function has been a
common topic in number theory for a long time. For the related background
information, we refer the readers to \cite{[A]} and \cite{[M]}.
Recently, this topic is still active. When studying the
arithmetic properties of the least common multiple of finitely
many consecutive positive integers,
Farhi \cite{[F]} defined the arithmetic function $g_k$ for any
positive integer $n$ by $g_k(n):=\frac{\prod_{i=0}^{k}(n+i)}{{\rm
lcm}_{0\le i\le k}\{n+i\}}$. In the same paper, Farhi showed that
$g_k$ is periodic of $k!$ and posed an open problem of determining
the smallest period of $g_k$. Let $P_k$ be the smallest period of
$g_k$. Define $L_0:=1$ and for any integer $k\ge 1$, we define
$L_k:={\rm lcm}(1,...,k)$. Subsequently, Hong and Yang \cite{[HY]}
improved the period $k!$ to $L_k$ and produced a conjecture stating that
$\frac{L_{k+1}}{k+1}$ divides $P_k$ for all nonnegative integers
$k$. By proving the Hong-Yang conjecture, Farhi and Kane \cite{[FK]}
determined the smallest period of $g_k$ and finally solved the open
problem posed by Farhi \cite{[F]}.
Let $k\ge 1,a\ge 1,b\ge 0$ and $ c\ge 1$ be integers. Let ${\mathbb
Q}$ and ${\mathbb N}$ denote the field of rational numbers and the
set of nonnegative integers. Define ${\mathbb N}^*:={\mathbb
N}\setminus\{0\}$. In order to investigate
the least common multiple of any $k+1$ consecutive terms in
the arithmetic progression $\{b+am)\}_{m\in \mathbb{N}^*}$,
Hong and Qian \cite{[HQ]} introduced the arithmetic function
$g_{k,a,b}$ defined for any positive integer $n$ by
$g_{k,a,b}(n):=\frac{\prod_{i=0}^k (b+a(n+i))}{{\rm lcm}_{0\le i\le
k}\{b+a(n+i)\}}.$ They \cite{[HQ]} showed that $g_{k,a,b}$ is periodic
and obtained the formula of the smallest period of $g_{k,a,b}$, which
extends the Farhi-Kane theorem to the general arithmetic progression case.

Let $f$ be a multiplicative function with $f(n)\ne 0$ for all
$n\in\mathbb{N}^*$. To measure the difference between
$\prod_{i=0}^kf(b+a(n+ic))$ and $f({\rm lcm}_{0\le i\le
k}\{b+a(n+ic)\})$, we define the arithmetic function $g_{k,f}$ for
any positive integer $n$ by
$$
g_{k,f}(n):=\frac{\prod_{i=0}^k f(b+a(n+ic))} {f({\rm lcm}_{0\le
i\le k} \{b+a(n+ic)\})}.       \eqno (1.1)
$$
One naturally asks the following interesting question.\\

\noindent{\bf Problem 1.1.} Let $f$ be a multiplicative function
such that $f(n)\ne 0$ for all positive integers $n$. Is $g_{k,f}$
periodic, and if so, what is the smallest period of $g_{k, f}$?\\

As usual, for any prime number $p$, we let $v_{p}$ be the normalized
$p$-adic valuation of ${\mathbb Q}$, i.e., $v_p(a)=s$ if
$p^{s}\parallel a$. For any real number $x$, by $\lfloor x\rfloor$
we denote the largest integer no more than $x$. Evidently,
$v_p(L_k)={\rm max}_{1\leq i\leq k}\{v_{p}(i)\}=\lfloor {\rm
log}_{p}k\rfloor$ is the largest exponent of a power of $p$ that is
at most $k$. We have the first main result of this paper which
answers the first part of Problem 1.1.\\

\noindent{\bf Theorem  1.2.} {\it Let $k\ge 1, a\ge 1, b\ge 0$ and
$c\ge 1$ be integers.  If $f$ is a multiplicative function so that
$f(n)\neq 0$ for all $n\in \mathbb{N}^*$, then the arithmetic
function $g_{k,f}$ is periodic and $cL_k$ is its period.}\\

It seems to be difficult to answer completely the second part of Problem 1.1.
We here are able to answer it for the Euler phi function $\varphi $ case. In fact,
we first prove a generalization of Hua's identity and then use it to
show that the arithmetic function $g_{k, \varphi}$ is periodic. Subsequently,
we develop $p$-adic techniques to determine the exact value of the smallest period
of $g_{k, \varphi}$. Note that it was proved by Farhi and Kane \cite{[FK]}
that there is at most one prime $p\le k$ such that $v_p(k+1)\ge v_p(L_k)\ge 1$.
We can now state the second main result of this paper as follows.\\

\noindent{\bf Theorem 1.3.} {\it Let $k\ge 1, a\ge 1, b\ge 0$ and
$c\ge 1$ be integers. Let $d:={\rm gcd}(a,b)$ and $a':=a/d$. Then
$g_{k,\varphi}$ is periodic, and its smallest period equals $Q_{k,
a', c}$ except that $v_p(k+1)\ge v_p(L_k)\ge 1$ for at most one odd
prime $p\nmid a'$, in which case its smallest period is equal to
$\frac{Q_{k, a', c}}{p^{v_p(L_k)}}$, where
$$Q_{k, a', c}:=\frac{cL_k}{\eta _{2,k,a',c}\prod_{{\rm prime}
\ q|a'}q^{v_{q}(cL_k)}},\eqno(1.2)$$ and
\begin{align*}
\eta_{2,k,a',c}:={\left\{
\begin{array}{rl}
2^{v_2(L_k)}, &\text{if} \ 2\nmid a' \ {\text and} \ v_2(k+1)\ge v_2(L_k)\ge 2, \\
2, &\text{if} \  2\nmid a \ and\ v_2(cL_k)=1, {\text or} \ k=3,
2\nmid a\ {\text and} \ 2|c,  {\text or} \ k=3, 2\nmid a'\ {\it
and}\ 2|d,\\
1, &\text{otherwise}.
\end{array}
\right.}
\end{align*}}
\\
So we answer the second part of Problem 1.1 for the Euler phi function.

The paper is organized as follows. In Section 2, we show that
$g_{k,f}$ is periodic and $cL_k$ is its period. In Section 3, we
provide a detailed $p$-adic analysis to the periodic arithmetic
function $g_{k,\varphi}$, and finally we determine the smallest
period of $g_{k,\varphi}$. The final section is devoted to the proof
of Theorem 1.3.

\section{\bf Proof of Theorem 1.2}

In this section, we give the proof of Theorem 1.2. We begin with
the following lemma.\\

\noindent{\bf Lemma 2.1.} {\it Let $A$ be any given totally ordered
set, and $a_1, ..., a_n$ be any  $n$ nonzero elements of $A$ (not
necessarily different). If we can define formal multiplication and
formal division for the set $A$, then we have
$$\max(a_1, ..., a_n)=a_1\cdots a_n \prod_{r=2}^n\prod_
{1\le i_1<\cdots<i_r\le
n}(\min(a_{i_1},\ldots,a_{i_r}))^{(-1)^{r-1}}.$$}

\begin{proof}
Rearrange these $n$ elements $a_1, ..., a_n$ such that $a_{j_1}\ge
\cdots\ge a_{j_n}$. For convenience, we let $b_i=a_{j_i}, \
i=1,2,\ldots,n$. Then the desired result in  Lemma 2.1  becomes
$$b_1=b_1\cdots b_n\prod_{r=2}^n\prod_{1\le i_1<\cdots<i_r\le n}
(\min(b_{i_1},\ldots,b_{i_r}))^{(-1)^{r-1}}.\eqno (2.1)$$ To prove
the result, it suffices to prove that for each $b_i$ the number of
times that $b_i$ occurs on the left side of (2.1) equals the number
of times that $b_i$ occurs on the right side of (2.1). We
distinguish the following two cases.

{\sc Case 1.} If $b_1=b_2=\cdots =b_n$, then the number of times that
$b_1$ occurs on the right side of (2.1) is
$$n-{n\choose 2}+\cdots+(-1)^{n-1}{n\choose n}=
-1+\sum_{r=1}^n(-1)^{r-1}{n\choose r}+1=-(1-1)^n+1=1.$$ Whereas, 1
is just the number of times that $b_1$ occurs on the left side of
(2.1).

{\sc Case 2.} If there exists a positive integer $s< n$ such that
$b_1=b_2=\cdots=b_s>b_{s+1}$, then the number of times $b_1$ occurs
on the right side of (2.1) is: $s-{s\choose
2}+\cdots+(-1)^{s-1}{s\choose s}=1$. For any $j>s$, we can always
assume that $b_{t+1},...,b_j,..., b_{t+l}$ are just the $l$ terms of
the sequence $\{b_i\}_{i=1}^n$ such that
$b_s>b_{t+1}=\cdots=b_j=\cdots=b_{t+l}$ for some $t\ge s$. Thus, the
number of times that $b_j$ occurs on the right side of (2.1) is
\begin{align*}
&l-\Big({t+l\choose 2}-{t\choose
2}\Big)+\cdots+(-1)^{t-1}\Big({t+l\choose t}-{t\choose
t}\Big)+(-1)^{t}{t+l\choose t+1}+\cdots+(-1)^{t+l-1}{t+l\choose t+l}\\
&=l+\sum_{r=2}^t(-1)^r{t\choose
r}+\sum_{i=2}^{t+l}(-1)^{i-1}{t+l\choose i}\\
&=l+(1-1)^t+{t\choose 1}-{t\choose 0}-(1-1)^{t+l}+{t+l\choose
0}-{t+l\choose 1}=0.
\end{align*}
This completes the proof of Lemma 2.1.
\end{proof}

In \cite{[Hu]}, Hua gave the following beautiful identity
$$ {\rm
lcm}(a_1, ..., a_n)=a_1\cdots a_n \prod_{r=2}^{n}\prod_{1\leq
i_1<\cdots<i_{r}\leq n}({\rm gcd}(a_{i_1}, ...,
a_{i_{r}}))^{(-1)^{r-1}},$$
where $a_1, ..., a_n$ are any given $n$
positive integers. In what follows, using Lemma 2.1, we generalize
the above Hua's identity to the multiplicative function case.\\

\noindent{\bf Lemma 2.2.} {\it Let $f$ be a multiplicative function,
and $a_1,a_2,\ldots,a_n$ be any $n$  positive integers. If $f(m)\ne
0$ for each $m\in \mathbb{N}^*$, then
$$f({\rm lcm}(a_1,a_2,\ldots,a_n))=f(a_1)\cdots f(a_n)\cdot
\prod_{r=2}^{n}\prod_{1\leq i_1<\cdots<i_{r}\leq
n}\big(f({\gcd}(a_{i_1}, ..., a_{i_{r}}))\big)^{(-1)^{r-1}}.$$ }

\begin{proof}
 Since
$f$ is a multiplicative function, we have
$$f({\rm lcm}(a_1,a_2,\ldots,a_n))=\prod_{p \ {\rm prime}}
f(p^{\max(v_p(a_1),v_p(a_2),\ldots,v_p(a_n))})$$ and
$$f(\gcd(a_{i_1},\ldots,a_{i_r}))=\prod_{p \ {\rm prime}}
f(p^{\min (v_p(a_{i_1}),\ldots,v_p(a_{i_r}))}).$$ Thus it suffices
to prove that
\begin{align*}
f(p^{\max(v_p(a_1),v_p(a_2),\ldots,v_p(a_n))})=
\prod_{r=1}^{n}\prod_{1\leq i_1<\cdots<i_{r}\leq
n}\big(f(p^{\min((v_p(a_{i_1}), ...,
v_p(a_{i_{r}}))})\big)^{(-1)^{r-1}}\  (2.2)
 \end{align*} for every
prime $p$. Now we define an order $\succeq$ for the set $S=\{f(p^m):
\ m\in \mathbb{N}\}$ according to the size of the power $m$ of the
prime $p$. That is, $f(p^i)\succeq f(p^j)$ if $i\ge j$ and
$f(p^i)\succ f(p^j)$ if $i> j$. It is easy to check that $S$ is a
totally ordered set for the order $\succeq$. So the equality (2.2)
follows immediately from Lemma 2.1 by letting $a_i$ be $f(p^{v_p(a_i)})$
for $1\le i\le n$ in Lemma 2.1. The proof of Lemma 2.2 is complete.
\end{proof}

If ${\rm gcd}(a_{i}, a_{j})={\rm gcd}(b_{i }, b_{j})$ for any $1\leq
i<j\le n$, then for any $t\ge 3$, one has ${\rm
gcd}(a_{i_1},a_{i_2},\ldots,a_{i_t})={\rm
gcd}(b_{i_1},b_{i_2},\ldots,b_{i_t})$ for any $1\leq
i_1<\cdots<i_t\leq n$. Therefore we immediately derive the following
result from Lemma 2.2.\\

\noindent{\bf Lemma 2.3.} {\it Let $a_1,a_2,\ldots,a_n$ and
$b_1,b_2,\ldots,b_n$ be any $2n$ positive integers. Let $f$ be a
multiplicative function with $f(n)\ne 0$ for all $n\in
\mathbb{N}^*$. If ${\rm gcd}(a_{i}, a_{j})={\rm gcd}(b_{i }, b_{j})$
for any $1\leq i<j\leq n$, then we have
\begin{align*}
\frac{\prod_{1\le i\le n}f(a_i)}{f({\rm lcm}_{1\le i\le n}\{a_i\})}
=\frac{\prod_{1\le i\le n}f(b_i)}{f({\rm lcm}_{1\le i\le
n}\{b_i\})}.
\end{align*}}\\

We are now in a position to show Theorem 1.2.\\
\\
{\it Proof of Theorem 1.2.} Let $n$ be a given  positive integer.
For any $0\leq i<j\leq k$, we have
\begin{align*}
{\rm gcd}(b+a(n+ic+cL_k),b+a(n+jc+cL_k))&={\rm gcd}(b+a(n+ic+cL_k),(j-i)ac)\\
&={\rm gcd}(b+a(n+ic),(j-i)ac)\\
&={\rm gcd}(b+a(n+ic),b+a(n+jc)).
\end{align*}
Thus by Lemma 2.3, we obtain that $g_{k,f}(n+cL_k)=g_{k,f}(n)$ for
any positive integer $n$. Therefore  $g_{k,f}$ is periodic and
$cL_k$ is
its period. \hfill $\square$\\

Obviously, by Theorem 1.2, the arithmetic function $g_{k,\varphi}$
is periodic and $cL_k$ is its period. In the next section, we will
provide detailed $p$-adic analysis to the arithmetic function $g_{k,
\varphi}$ which leads us to determine the exact value of the
smallest period of $g_{k, \varphi}$.

\section{\bf Local analysis of $g_{k,\varphi}$}

Throughout this section, we let $a'=a/d$ and $ b'=b/d$ with
$d=\gcd(a, b)$. Then $\gcd(a',b')=1$. Let
$$
S_{k,a',b',c}(n):=\{b'+a'n,b'+a'(n+c),\ldots,b'+a'(n+kc)\}
$$
be the set consisting of  $k+1$ consecutive jumping terms with gap
$c$ in the arithmetic progression $\{b'+a'm\}_{m\in \mathbb{N}}$.
For any given prime number $p$, define $g_{p,k,\varphi}$ for any
$n\in \mathbb {N}^*$ by
$g_{p,k,\varphi}(n):=v_{p}(g_{k,\varphi}(n))$. Let $P_{k, \varphi}$
be the smallest period of $g_{k,\varphi}$. Then $g_{p,k,\varphi}$ is
a periodic function for each prime $p$ and $P_{k,\varphi}$ is a
period of $g_{p,k,\varphi}$. Let $P_{p, k, \varphi}$ be the smallest
period of $g_{p,k,\varphi}$. Since
\begin{align*}
\varphi(b+a(n+ic))&=\varphi(d(b'+a'(n+ic)))=\varphi\big(
\prod_{p|d}p^{v_p(b'+a'(n+ic))+v_p(d)}\prod_{p\nmid
d}p^{v_p(b'+a'(n+ic))}\big)\\
&=\bigg(\prod_{p|d}p^{v_p(d)-1}(p-1)\bigg)\bigg
(\prod_{p|d}p^{v_p(b'+a'(n+ic))}\bigg)\bigg(\prod_{p\nmid
d}\varphi\big(p^{v_p(b'+a'(n+ic))}\big)\bigg)\\
&=\varphi(d)\bigg(\prod_{p|d}p^{v_p(b'+a'(n+ic))}\bigg)\bigg(\prod_{p\nmid
d}\varphi\big(p^{v_p(b'+a'(n+ic))}\big)\bigg),
\end{align*}
we have
\begin{align*}
 g_{k,\varphi}(n)&=\frac{\prod_{i=0}^k \varphi(b+a(n+ic))} {\varphi({\rm lcm}_{0\le
i\le k} \{b+a(n+ic)\})}=\frac{\prod_{i=0}^k \varphi(d(b'+a'(n+ic)))}
{\varphi(d\cdot {\rm lcm}_{0\le
i\le k} \{b'+a'(n+ic)\})}\\
&= \frac{\prod_{i=0}^k\Big(\varphi(d)\Big(\prod_{p| d}
p^{v_p(b'+a'(n+ic))}\Big)\Big(\prod_{p\nmid
d}\varphi(p^{v_p(b'+a'(n+ic))})\Big)\Big)}{\varphi(d)\Big(\prod_{p|
d}p^{\max_{0\le i\le k}\{v_p(b'+a'(n+ic))\}}\Big)\Big(\prod_{p\nmid
d}\varphi(p^{\max_{0\le i\le k}\{v_p(b'+a'(n+ic))\}})\Big)}\\
&=(\varphi(d))^k\frac{\prod_{i=0}^k \Big(\Big(\prod_{p| d}
p^{v_p(b'+a'(n+ic))}\Big)\Big(\prod_{p\nmid
d}\varphi(p^{v_p(b'+a'(n+ic))})\Big)\Big)}{ \Big(\prod_{p|
d}p^{\max_{0\le i\le k}\{v_p(b'+a'(n+ic))\}}\Big)\Big(\prod_{p\nmid
d} \varphi(p^{\max_{0\le i\le k}\{v_p(b'+a'(n+ic))\}})\Big)}.
\end{align*}
Note that for any prime $q$, we have that for any positive integer
$e$, $\varphi(q^e)=q^{e-1}(q-1)$. So when computing  $p$-adic
valuation of  $g_{k, \varphi}(n)$, we not only need to compute
$v_p(\varphi (p^{\alpha }))$ for $\alpha \ge 2$, but also need to
consider $p$-adic valuation of $q-1$ for those primes $q$ with
$p|(q-1)$. By some computations, we obtain the following two
equalities.

If $p\nmid d$, then
\begin{align*}
g_{p,k,\varphi}(n)&= \sum_{ m\in S_{k,a',b',c}(n)}\max(v_p(m)-1,0)-
\max(\max_{ m\in
S_{k,a',b',c}(n)}\{v_p(m)-1\},0)\\
& +\sum_{{\rm prime}\ q:\ q\nmid d, q\neq p}\max(0,\#\{m\in S_{k,a',b',c}(n): q| m\}-1)\cdot v_p(q-1)+kv_p(\varphi(d))\\
&= kv_p(\varphi(d))+\sum_{e\ge 2}{\max} (0, \#\{m\in
S_{k,a',b',c}(n):p^e|
m\}-1)\\
& \  \ +\sum_{{\rm prime} \ q:\ q\nmid d, p|(q-1)}\max(0,\#\{m\in
S_{k,a',b',c}(n): q| m\}-1)\cdot v_p(q-1). \  \  \  \  \  \  \
 \ (3.1)
\end{align*}

If $p| d$, then
\begin{align*}
g_{p,k,\varphi}(n)&=
kv_p(\varphi(d))+\sum_{i=0}^kv_p(b'+a'(n+ic))-\max_{0\le i\le
k}\{v_p(b'+a'(n+ic))\}\\
& \  \  +\sum_{{\rm prime}\ q:\ q\nmid d, q\neq p}\max(0,\#\{m\in S_{k,a',b',c}(n): q| m\}-1)\cdot v_p(q-1)\\
&=kv_p(\varphi(d))+\sum_{e\ge 1}{\max} (0, \#\{m\in
S_{k,a',b',c}(n):p^e|
m\}-1)\\
& \  \ +\sum_{{\rm prime} \ q:\ q\nmid d, p|(q-1)}\max(0,\#\{m\in
S_{k,a',b',c}(n): q| m\}-1)\cdot v_p(q-1). \ \ \  \ \ \
 \ (3.2)
\end{align*}

In order to analyze the function $g_{p,k,\varphi}$ in detail, we
need the following results.\\

\noindent{\bf Lemma 3.1.} {\it  Let $e$ and $m$ be positive
integers. If $p\nmid a'$, then any $p^e$ consecutive terms in the
arithmetic progression $\{b'+a'(m+ic)\}_{i\in \mathbb{N}}$ are
pairwise incongruent modulo $p^{v_p(c)+e}$. In particular, there is
at most one term  divisible by $p^e$ in $S_{k,a',b',c}(n)$ for
$e>v_p(cL_k)$.}
\begin{proof}
Suppose that there are two integers $i,j$ such that $0<j-i\le p^e-1$
and $b'+(m+ic)a'\equiv b'+(m+jc)a' \pmod {p^{v_p(c)+e}}$.
Then $p^e\mid (j-i)a'$. Since ${\rm
gcd}(p,a')=1$, we have $p^{e}\mid (j-i)$. This is a contradiction.
\end{proof}

\noindent{\bf Lemma 3.2.} {\it Let $F$ be a positive rational-valued
arithmetic function. For any prime $p$, define $F_p$ by
$F_p(n):=v_p(F(n))$ for $n\in {\mathbb N}^*$. Then $F$ is periodic
if and only if $F_p$ is periodic for each prime $p$ and ${\rm
lcm}_{{\rm prime} \ p}\{T_{p,F}\}$ is finite, where $T_{p,F}$ is the
smallest period of $F_p$. Furthermore, if $F$ is periodic, then the
smallest period $T_F$ of $F$ is equal to ${\rm lcm}_{{\rm prime} \
p}\{T_{p,F}\}$.}
\begin{proof}
$\Rightarrow)$ Since $F$ is periodic and $T_F$ is its smallest
period, we have $F(n+T_F)=F(n)$ for any $n\in \mathbb{N}^*$, and
hence $F_p(n+T_F)=v_{p}(F(n+T_F))=v_{p}(F(n))=F_p(n)$. In other
words, $F_p$ is periodic and $T_F$ is a period of $F_p$ for every
prime $p$. So we have ${\rm lcm}_{{\rm prime} \ p}\{T_{p,F}\}|
T_{F}$ and ${\rm lcm}_{{\rm prime} \ p}\{T_{p,F}\}$ is finite.

$\Leftarrow)$ Since for any $n\in \mathbb{N}^*$, we have that
$v_{q}(F(n+{\rm lcm}_{{\rm prime} \ p \ }\{T_{p,F}\}))=v_{q}(F(n))$
for each prime $q$. Thus  $F(n+{\rm lcm}_{{\rm prime} \
p}\{T_{p,F}\})=F(n)$ for any $n\in \mathbb{N}^*$. So $F$ is periodic
and ${\rm lcm}_{{\rm prime} \ p}\{T_{p,F}\}$ is a period of it.
Hence  $T_{F}$ divides ${\rm lcm}_{{\rm prime} \ p}\{T_{p,F}\}$.

From the above discussion, we immediately derive that $T_{F}={\rm
lcm}_{{\rm prime} \ p}\{T_{p,F}\}$ if $F$ is periodic.
\end{proof}

For any prime $p\ge cL_k+1$, we have by Lemma 3.1 that there is at
most one term divisible by $p$ in $S_{k,a',b',c}(n)$ and there is at
most one element divisible by the  prime $q$ satisfying $p| (q-1)$
in $S_{k,a',b',c}(n)$. Thus for any prime $p\ge cL_k+1$, we can get
from (3.1) and (3.2) that $g_{p,k,\varphi}(n)=kv_p(\varphi(d))$ for
every positive integer $n$. Namely, we have $P_{p,k,\varphi}=1$ for
each prime $p$ such that $p\ge cL_k+1$. Thus  by Lemma 3.2, we
immediately have the following.\\

\noindent{\bf Lemma 3.3.} {\it We have}
$$
P_{k,\varphi}={\rm lcm}_{{\rm prime} \ p\le cL_k}
\{P_{p,k,\varphi}\}.
$$\\

In what follows it is enough to compute $P_{p,k,\varphi}$ for every
prime $p$ with $p\le cL_k$.  First we need to simplify
$g_{p,k,\varphi}$ for  $p\le cL_k$. For any prime $q$ satisfying
$q\nmid cL_k$, we obtain by Lemma 3.1 that there is at most one term
divisible by $q$ in $S_{k,a',b',c}(n)$. On the other hand, for any
prime $q$ satisfying $q| a'$,  we have $\gcd (q,b')=1$ since ${\rm
gcd}(a',b')=1$. Thus for $0\leq i\leq k$, we have that $\gcd(q,
b'+a'(n+ic))=1$ for all  $n\in \mathbb{N}^*$. So there is no term
divisible by any prime factor $q$ of  $a'$ in
$S_{k,a',b',c}(n)$. Thus from (3.1) and (3.2), we derive the following equality:\\
$$g_{p,k,\varphi}(n)=kv_p(\varphi(d))+\sum_{e=1}^{v_p(cL_k)}f_{e}(n)+
\sum_{{\rm prime} \ q: \ q| cL_k \atop p|(q-1), \ q\nmid a}h_q(n),
\eqno (3.3)
$$
where
\begin{align*}
f_{e}(n):={\left\{
  \begin{array}{rl}
0, \quad&\text{if} \ p\nmid d\ \text{and}\ e= 1,\\
\max (0, \#\{m\in S_{k,a',b',c}(n):p^e| m\}-1),
\quad&\text{otherwise}
 \end{array}
\right.}
\end{align*}
and $$ h_q(n):=\max(0,\#\{m\in S_{k,a',b',c}(n): q| m\}-1)\cdot
v_p(q-1).$$ For any positive integer $n$, it is easy to check that
$f_{e}(n+p^{v_p(cL_k)})=f_{e}(n)$ for each $1\le e\le v_p(cL_k)$ and
$h_q(n+q)=h_q(n)$ for each prime $q$ such that $q\nmid a, \ q|
cL_k$ and $p| (q-1)$. Consequently, we obtain that
$p^{v_p(cL_k)}\prod_{{\rm prime}\ q:\  q| cL_k \atop q\nmid a, \ p|
(q-1)}q$ is a period of the function $g_{p,k,\varphi}$.  To get the
smallest period of $g_{p,k,\varphi}$ for each prime $p\le cL_k$, we
need to make  more detailed $p$-adic analysis about
$g_{p,k,\varphi}$.
We divide it into the following four cases.\\

\noindent{\bf Lemma 3.4.} {\it Let $p$ be a prime such that $p\le
cL_k$ and  $p\nmid cL_k$. Then
\begin{align*}
P_{p,k,\varphi}= \prod _{{\rm prime} \ q :\ q| c \atop q\nmid a, \
p| (q-1)} q.
\end{align*}
 }
\begin{proof}
Since $p\le cL_k$ and  $p\nmid cL_k$, we have  $k<p\le cL_k$ and
$v_p(cL_k)=0$. Hence we have that
$g_{p,k,\varphi}(n)=kv_p(\varphi(d))+\sum_{{\rm prime}\ q: \ q| cL_k
\atop p| (q-1), \ q\nmid a}h_q(n)$
by (3.3). If there is no prime
$q$ satisfying $q| cL_k$, $p| (q-1)$ and $q\nmid a$, then we have
$g_{p,k,\varphi}(n)=kv_p(\varphi(d))$ for any positive integer $n$,
and hence $P_{p,k,\varphi}=1$ for such primes $p$. If there is a
prime $q$ satisfying $q| cL_k$, $p| (q-1)$ and $q\nmid a$, then we
must have $q| c$, since such $q$ must satisfy $q\ge p+1>k$. So by
the argument before Lemma 3.4, we have that $A:=\prod_{ \ {\rm prime} \
q : \ q| c\atop q\nmid a, \ p| (q-1)}q$ is a period of
$g_{p,k,\varphi}$.

Now it remains to prove that  $A$ is the smallest
period of $g_{p,k,\varphi}$. For any prime factor $q$ of $A$, we can
choose a positive integer $n_0$ such that $v_{q}(b'+a'n_0)\ge 1$
because $q\nmid a$. Since $q>k$ and $q| c$, we have that
$v_q(b'+a'(n_0+ic))\ge 1$ and $v_q(b'+a'(n_0+ic+A/q))=v_q(a'A/q)=0$ for each $0\le i\le k$.
Hence  there is no term divisible by $q$ in
$S_{k,a',b',c}(n_0+A/q)$. Thus  $h_q(n_0)=kv_p(q-1)\ge
k>0=h_q(n_0+A/q)$. On the other hand,  $h_{q'}(n_0)=h_{q'}(n_0+A/q)$
for any other prime factors $q'\ne q$ of $A$. It then follows that
$g_{p,k,\varphi}(n_0)\ne g_{p,k,\varphi}(n_0+A/q)$.  Therefore  $A$
is the smallest period of $g_{p,k,\varphi}$. This completes the
proof of Lemma 3.4.
\end{proof}

\noindent{\bf Lemma 3.5.} {\it Let $p$ be a prime such that $p|
cL_k$ and $p|a'$. Then
\begin{align*}
P_{p,k,\varphi}=\bigg(\prod_{{\rm prime} \ q:\  q\nmid ac, \ q| L_k
\atop p| (q-1), \ k+1\not\equiv 0\pmod {q}} \ q
\bigg)\bigg(\prod_{{\rm prime}\ q:\ q\nmid a,\atop q| c,\
p|(q-1)}q\bigg).
\end{align*}
}
\begin{proof}
For convenience, we let $A$ denote the number on the right side of
the equality in Lemma 3.5. First, we prove that $A$ is a period of
$g_{p,k,\varphi}$. Since there is no term divisible by $p$ in
$S_{k,a',b',c}(n)$ if $p| a'$, we have that $f_e(n)=0$ for any
positive integer $n$ and for each $1\le e\le v_p(cL_k)$. Hence we
have $g_{p,k,\varphi}(n)=kv_p(\varphi(d))+\sum_{{\rm prime}\ q:\ q|
cL_k \atop p| (q-1), \ q\nmid a}h_q(n)$. If there is no prime $q$
such that $q| cL_k$, $p| (q-1)$ and $q\nmid a$, then
$g_{p,k,\varphi}(n)=kv_p(\varphi(d))$ for any positive integer $n$,
and hence $P_{p,k,\varphi}=1$ for such primes $p$.  If $q\nmid ac,
q|L_k$ and $k+1\equiv0\pmod{q}$, then  any $q$ consecutive terms in
the arithmetic progression $\{b'+a'(m+ic)\}_{i\in \mathbb{N}}$ are
pairwise incongruent modulo $q$ by Lemma 3.1. Therefore, we have
$h_q(n)=h_q(n+1)=v_p(q-1)(\frac{k+1}{q}-1)$ for any positive integer
$n$. Namely, 1 is a period of $h_q$ for such primes $q$. Since $q$
is a period of $h_q$ for any other  primes $q$ such that  $q| cL_k,
q\nmid a$ and $ p| (q-1)$, we have that $A$ is a
period of the function $g_{p,k,\varphi}$. To prove that $A$ is the
smallest period of $g_{p,k,\varphi}$, it is enough to show that
$A/q$ is not the period of $g_{p,k,\varphi}$ for every prime factor
$q$ of $A$. We divide the prime factors of $A$ into the following
two cases.

{\sc Case 1.} $q$ is a prime factor  of $A$ such that $q\nmid ac$,
 $q| L_k$, $p| (q-1)$ and $k+1\not\equiv0\pmod{q}$. To
prove $A/q$ is not the period of $g_{p,k,\varphi}$, it suffices to
prove that $A/q$ is not the period of the function $h_q$ since $A/q$
is a period of $h_{q'}$ for any other primes $q'\ne q$ such that $q'|
cL_k, p| (q'-1)$ and $q'\nmid a$. Since $\gcd(A/q, q)=1$, there
exists a positive integer $r_0$ such that $r_0A/q\equiv 1\pmod{q}$.
We pick a positive integer $n_0$ so that $v_{q}(b'+a'n_0)\ge 1$
since $q\nmid a$. So we have that the terms divisible by $q$ in the
arithmetic progression $\{b'+a'(n_0+ic)\}_{i\in \mathbb{N}}$ must be
of the form  $b'+a'(n_0+tcq)$ for some $t\in \mathbb{N}$, and there
are at least two terms  divisible by $q$ in $S_{k,a',b',c}(n_0)$
since $q|L_k$. Comparing $S_{k,a',b',c}(n_0)$ with
$S_{k,a',b',c}(n_0+r_0A/q)$, we obtain that $b'+a'(n_0+j)\equiv
b'+a'(n_0+r_0A/q+j-1) \pmod{q} \ {\rm for \ each} \ 1\le j\le k$
while  $b'+a'(n_0+r_0A/q+k)\equiv b'+a'(n_0+k+1)\not\equiv
b'+a'n_0\equiv 0\pmod{q}$. Thus the number of terms divisible by $q$
in $S_{k,a',b',c}(n_0+r_0A/q)$ equals the number of terms divisible
by $q$ in $S_{k,a',b',c}(n_0)$ minus one, which means that
$h_q(n_0+r_0A/q)= h_q(n_0)-v_p(q-1)$. Therefore, $A/q$ is not the
period of $g_{p,k,\varphi}$ in this case.

{\sc Case 2.} $q$ is a prime factor  of $A$ satisfying $q\nmid a$,
$q| c$, and $p| (q-1)$. As above, we select two positive integers
$r_0$ and $n_0$ such that $r_0A/q\equiv1\pmod{q}$ and
$v_{q}(b'+a'n_0)\ge 1$. So we obtain that $v_q(b'+a'(n_0+ic))\ge 1$
and $v_q(b'+a'(n_0+r_0A/q+ic))=v_q(a'r_0A/q)=0$ for each $0\le i\le
k$. In other words, all the $k+1$ terms are divisible by $q$ in
$S_{k,a',b',c}(n_0)$ while no term is divisible by $q$ in
$S_{k,a',b',c}(n_0+r_0A/q)$.  Thus $h_q(n_0)=kv_p(q-1)\ge
k>0=h_q(n_0+r_0A/q)$. It follows immediately that $A/q$ is not the
period of $g_{p,k,\varphi}$ in this case.

So $A$ is the smallest period of $g_{p,k,\varphi}$. The proof of
Lemma 3.5 is complete.
\end{proof}

\noindent{\bf Lemma 3.6.} {\it Let $p$ be a prime such that $p|
cL_k$, $p\nmid a'$ and $p\nmid d$. Then
\begin{align*}
P_{p,k,\varphi}= p^{e(p,k)}\bigg(\prod_{{\rm prime} \ q:\ q\nmid
ac, \ q| L_k \atop p| (q-1), \ k+1\not\equiv 0\pmod {q}} \ q
\bigg)\bigg( \prod_{{\rm prime}\ q: \ q\nmid a,\atop q| c, \ p|
(q-1)}q\bigg),
\end{align*}
where
\begin{align*}
e(p,k):={\left\{
  \begin{array}{rl}
0, \quad&\text{if} \  v_p(cL_k)=1,\\
v_p(c), \quad&\text{if} \ v_{p}(k+1)\geq v_p(L_k) \ and \
v_p(cL_k)\ge 2,\\
v_p(cL_k), \quad&\text{if} \ v_p(k+1)< v_p(L_k) \ and \ v_p(cL_k)\ge
2.
 \end{array}
\right.}
\end{align*}
}
\begin{proof}
From (3.3), we get that
$$g_{p,k,\varphi}(n)=kv_p(\varphi(d))+\sum_{e= 2}^{v_p(cL_k)}f_e(n)+
\sum_{{\rm prime} \ q: \ q| cL_k \atop p| (q-1), \ q\nmid a}h_q(n).
\eqno (3.4)$$ Let $A$ denote the number $p^{e(p,k)}\big(\prod_{{\rm
prime} \ q:\ q\nmid ac, \ q| L_k \atop p| (q-1), \ k+1\not\equiv
0\pmod {q}} q\big)\big(\prod_{{\rm prime}\ q: \ q\nmid a\atop q| c, \ p|
(q-1)}q\big)$. We distinguish   the following three cases.

{\sc Case 1.} $v_p(cL_k)= 1$. Since $v_p(cL_k)= 1$, we have
$$
g_{p,k,\varphi}(n)=kv_p(\varphi(d))+\sum_{{\rm prime}\ q: \ q| cL_k
\atop p| (q-1), \ q\nmid a}h_q(n)
$$
for  every positive integer $n$ by (3.4). The process of proving
that $A$ is the smallest period of $g_{p,k,\varphi}$  is the same as
the proof of Lemma 3.5, one can easily check it.

{\sc Case 2.} $v_p(k+1)\ge v_p(L_k)$ and $v_p(cL_k)\ge 2$. We
consider the following two subcases.

{\sc Subcase 2.1.} $p>k$, $v_p(k+1)\ge v_p(L_k)$ and $ v_p(cL_k)\ge
2$. In this case, we have $v_p(cL_k)=v_p(c)\ge 2$. So we obtain that
$p^{v_p(c)}$ is a period of $f_e$ for each $2\le e\le
v_p(cL_k)=v_p(c)$. By the same method as in the proof of Lemma 3.5,
we can derive that $A/p^{v_p(c)}$ is a period of $\sum_{{\rm prime}
\ q:\ q| cL_k \atop p| (q-1), \ q\nmid a}h_q(n)$, and hence by
(3.4) $A$ is a period of $g_{p,k,\varphi}$. Now it suffices to
prove that  $A/P$ is not the period of $g_{p,k,\varphi}$ for any
prime factor $P$ of $A$.

For the prime $p$,  we have by (3.4) that $A/p$ is a period of $f_e$
for each $2\le e\le v_p(c)-1$  and is also a period of $h_q$ for
each prime $q$ such that $q\nmid a, \ q| cL_k$ and $ p| (q-1)$. So
it is enough to prove that $A/p$ is not the period of $f_{v_p(c)}$.
 Since $p\nmid a'$, we can choose a positive integer $n_0$ such that
$v_p(b'+a'n_0)=v_p(c).$ It is easy to see that $p^{v_p(c)}|
b'+a'(n_0+ic)$  and  $p^{v_p(c)}\nmid b'+a'(n_0+A/p+ic)$ for each
$0\le i\le k$. Thus comparing the two sets $S_{k,a',b',c}(n_0)=\{
b'+a'(n_{0}+ic)\}_{0\le i\le k}$ and
$S_{k,a',b',c}(n_0+A/p)=\{b'+a'(n_0+A/p+ic)\}_{0\le i\le k}$, we
obtain that
$$f_{v_p(c)}(n_0)= \max (0, \#\{m\in S_{k,f,c}(n_0):p^{v_p(c)}| m\}-1)=k,$$
\begin{align*}
f_{v_p(c)}(n_0+A/p)&= \max (0, \#\{m\in
S_{k,f,c}(n_0+A/p):p^{v_p(c)}| m\}-1)=0.
\end{align*}
Therefore, $A/p$ is not the period of $g_{p,k,\varphi}$.

For any prime factor $q$ of $A$ such that $q| cL_k$, $p| (q-1)$ and
$q\nmid a$. It is easy to see that $A/q$ is a period of $\sum_{e=
2}^{v_p(cL_k)}f_e(n)$ and  $ h_{q'}(n)$ for each prime factor $q'\ne
q$ of $A$ satisfying $ q'|cL_k, p|(q'-1)$ and $q'\nmid a$. Similarly
to the proof of Lemma 3.5, we can deduce that $A/q$ is not the
period of $h_q$, and hence $A/q$ is not the period of
$g_{p,k,\varphi}$. Therefore, $A$ is the smallest period of
$g_{p,k,\varphi}$ in this subcase.

{\sc Subcase 2.2.} $p\le k$, $v_p(k+1)\ge v_p(L_k)$ and $
v_p(cL_k)\ge 2$. To prove that  $A$ is a period of
$g_{p,k,\varphi}$, it suffices to prove that $p^{v_p(c)}$ is a
period of $f_e$ for each $1\le e\le v_p(cL_k)$ by the argument in
Subcase 2.1 of this proof.  For any given positive integer $n$,
comparing the two sets $S_{k,a',b',c}(n)=\{ b'+a'(n+ic)\}_{0\le i\le
k}$ and $S_{k,a',b',c}(n+c)=\{ b'+a'(n+ic)\}_{1\le i\le k+1}$, we
find that their distinct terms are $b'+a'n$ and $b'+a'(n+(k+1)c)$.
From $v_p(k+1)\ge v_p(L_k)$ we deduce that
$$b'+a'n\equiv b'+a'(n+(k+1)c)\pmod {p^{v_p(cL_k)}}.$$
Therefore we obtain that $f_e(n)=f_e(n+c)$ for each $e\in
\{2,\ldots, v_p(cL_k)\}$.  Since $\gcd(c/p^{v_p(c)},p)=1$, we
can  always find two  integers $t,t_1$ such that
$tc/p^{v_p(c)}=t_1{p^{v_p(L_k)}}+1$. Note that $p^{v_p(cL_k)}$ is a
period of $f_e$ for each $2\le e\le v_p(cL_k)$. Therefore, we have
$f_e(n+p^{v_p(c)})=f_e(n
+p^{v_p(c)}+t_1p^{v_p(cL_k)})=f_e(n+tcp^{v_p(c)}/p^{v_p(c)})=f_e(n+tc)=f_e(n)$
for each positive integer $n$ and each $2\le e\le v_p(cL_k)$.
 Thus $A$ is
a period of $g_{p,k,\varphi}$ as required.

 Now we only need to
prove that $A/P$ is not the period of $g_{p,k,\varphi}$ for any
prime factor $P$ of $A$. For any prime factor $q$ of $A$ such that
$q| cL_k$, $p| (q-1)$ and $q\nmid a$, the proof is similar to Subcase 2.1.
If $v_p(c)=0$, then $p$ is not a prime factor of
$A$, and the proof of this case is complete.   In the following, we
need to prove that $A/p$ is not the period of $g_{p,k,\varphi}$ if
$v_p(c)\ge 1$. Since $A/p$ is a period of $h_q$ for each $q$ such
that $q\nmid a, \ q| cL_k$ and $ p| (q-1)$, it is enough to prove
that $A/p$ is not the period of the function
$\sum_{e=2}^{v_p(cL_k)}f_{e}(n)$.

If $v_p(c)\ge 2$, we choose $n_0\in \mathbb{N}^*$ such that
$v_p(b'+a'n_0)=v_p(c)$. Comparing $S_{k,a',b',c}(n_0)=\{
b'+a'(n_{0}+ic)\}_{0\le i\le k}$ with $
S_{k,a',b',c}(n_0+A/p)=\{b'+a'(n_0+A/p+ic)\}_{0\le i\le k}$, we
obtain that each term of $S_{k,a',b',c}(n_0)$ is divisible by
$p^{v_p(c)}$, while there is no term divisible by $p^{v_p(c)}$ in
$S_{k,a',b',c}(n_0+A/p)$ since $v_p(b'+a'(n_0+A/p+ic))=\min\big(
v_p(b'+a'(n_0+ic)), v_p(a'A/p)\big)=v_p(c)-1$ for each $0\le i\le
k$.  So

$$\sum_{e=2}^{v_p(cL_k)}f_{e}(n_0)\ge k(v_p(c)-1)>k(v_p(c)-2)
=\sum_{e=2}^{v_p(cL_k)}f_{e}(n_0+A/p).$$

If $v_p(c)=1$, then $v_p(A/p)=0$.  Choosing $n_0\in \mathbb{N}^*$
such that $v_p(b'+a'n_0)=v_p(cL_k)$, we have that there is at least
two terms $b'+a'n_0$ and $b'+a'(n_0+p^{v_p(L_k)}c)$ divisible by
$p^{v_p(cL_k)}$ in $S_{k,a',b',c}(n_0)$ but no term is divisible by
$p$ in $S_{k,a',b',c}(n_0+A/p)$ since $v_p(b'+a'(n_0+A/p+ic))=0$ for
all $0\le i\le k$. Therefore, we have
$$\sum_{e=2}^{v_p(cL_k)}f_{e}(n_0)\ge v_p(cL_k)-1>0
=\sum_{e=2}^{v_p(cL_k)}f_{e}(n_0+A/p).$$ Thus $A/p$ is not the
period of $g_{p,k,\varphi}$ in this case.

{\sc Case 3.} $p\le k$, $v_p(k+1)< v_p(L_k)$ and $v_p(cL_k)\ge 2$.
By the discussion before Lemma 3.4, it is easy to get that  $A$ is a
period of $g_{p,k,\varphi}$. As above,  it suffices to prove that
$A/P$ is not the period of $g_{p,k,\varphi}$ for any prime factor
$P$ of $A$ in the following. By a similar argument as in
Subcase 2.1, we now only need to show that $A/p$ is not the period
of $f_{v_p(cL_k)}$. Since $v_p(A/p)=v_p(cL_k)-1$, we can select a
positive integer $r_0$ such that $r_0A/p\equiv
p^{v_p(cL_k)-1}\pmod{p^{v_p(cL_k)}}$. In the following, we  prove
that $p^{v_p(cL_k)-1}$ is not the period of $f_{v_p(cL_k)}$, from
which we can deduce that $A/p$ is not the period of
$g_{p,k,\varphi}$. Since $v_{p}(k+1)<v_p(L_k)$, we can always
suppose that $k+1\equiv r\pmod {p^{v_p(L_k)}} \ {\rm for \ some}\
1\leq r\leq p^{v_p(L_k)}-1.$  We distinguish the following two
subcases.

{\sc Subcase 3.1.} $1\le r\le p^{v_p(L_k)}-p^{v_p(L_k)-1}$. Choose a
positive integer $n_0$ such that $v_p(b'+a'n_0)\ge v_p(cL_k)$.
Compare the number of terms divisible by $p^{v_p(cL_k)}$ in the two
sets $S_{k,a',b',c}(n_0)=\{b'+a'(n_0+kc)\}_{0\le i\le k}$ and $
S_{k,a',b',c}(n_0+p^{v_p(L_k)-1}c)=\{b+a(n_0+(p^{v_p(L_k)-1}+i)c)\}_{0\le
i \le k}$.  Since $\{b'+a'(n_0+p^{v_p(L_k)-1}c),\ldots,
b'+a'(n_0+kc)\}$ is the intersection of $S_{k,a',b',c}(n_0)$ and
$S_{k,a',b',c}(n_0+p^{v_p(L_k)-1}c)$,  it suffices to compare the
set $\{b'+a'n_0,\ldots,b'+a'(n_0+(p^{v_p(L_k)-1}-1)c)\}$ with the
set $\{b'+a'(n_{0}+(k+1)c),\ldots,b'+a'(n_0+(k+p^{v_p(L_k)-1})c)\}$.
By Lemma 3.1, we know that the terms divisible by $p^{v_p(cL_k)}$ in
the arithmetic progression $\{b'+a'(n_0+ic)\}_{i\in \mathbb{N}}$ are
of the form $b'+a'(n_0+tp^{v_p(L_k)}c), \ t\in \mathbb{N}$. Since
$k+1\equiv r\pmod {p^{v_p(L_k)}}$ and $1\leq r\leq
p^{v_p(L_k)}-p^{v_p(L_k)-1}$, we have $k+j\equiv r+j-1\not\equiv 0
\pmod {p^{v_p(L_k)}}$ for all $1\leq j\leq p^{v_p(L_k)-1}$. Hence
$p^{v_p(cL_k)}\nmid (b'+a'(n_0+(k+j)c))$ for all $1\leq j\leq
p^{v_p(L_k)-1}$. Whereas, $b'+a'n_0$ is the only term in the set
$\{b'+a'n_0,b'+a'(n_{0}+c),\ldots,b'+a'(n_{0}+(p^{v_p(L_k)-1}-1)c)\}$
which is divisible by $p^{v_p(cL_k)}$. Therefore we have
$$
f_{v_p(cL_k)}(n_0+p^{v_p(L_k)-1}c)=f_{v_p(cL_k)}(n_0)-1.$$

{\sc Case 3.2.} $p^{v_p(L_k)}-p^{v_p(L_k)-1}<r\leq p^{v_p(L_k)}-1$.
Pick a positive integer $n_0$ such that
$v_p(b'+a'(n_0+(p^{v_p(L_k)-1}-1)c))\ge v_p(cL_k)$. Then the terms
divisible by $p^{v_p(cL_k)}$ in the arithmetic progression
$\{b'+a'(n_0+ic)\}_{i\in \mathbb{N}}$ should be of the form
$b'+a'(n_0+(p^{v_p(L_k)-1}-1+tp^{v_p(L_k)})c)$, where $t\in
\mathbb{N}$. As in the discussion of Case 3.1, it is sufficient to
compare $\{b'+a'n_0,\ldots,b'+a'(n_0+(p^{v_p(L_k)-1}-1)c)\}$ with
$\{b'+a'(n_0+(k+1)c),\ldots,b'+a'(n_0+(k+p^{v_p(L_k)-1})c)\}$.  By
comparison, we obtain  that $p^{v_p(cL_k)}\nmid (b'+a'(n_0+(k+j)c))$
for all $1\leq j\leq p^{v_p(L_k)-1}$, while the term
$b'+a'(n_0+(p^{v_p(L_k)-1}-1)c)$ is the only term divisible by
$p^{v_p(cL_k)}$ in the set
$\{b'+a'n_0,\ldots,b'+a'(n_0+(p^{v_p(L_k)-1}-1)c)\}$. Hence we have
$f_{v_p(cL_k)}(n_0+p^{v_p(L_k)-1}c)=f_{v_p(cL_k)}(n_0)-1$. From the
argument in the above two subcases, we deduce that $p^{v_p(L_k)-1}c$
is not the period of $f_{v_p(cL_k)}$, which shows that $A/p$ is not
the period of $g_{p,k,\varphi}$ in Case 3.

 Thus $A$ is the smallest
period of $g_{p,k,\varphi}$  as desired. This completes the proof of
Lemma 3.6.
\end{proof}

\noindent{\bf Lemma 3.7.} {\it Let $p$ be a prime such that $p|
cL_k$, $p\nmid a'$  and $p| d$. Then
\begin{align*}
P_{p,k,\varphi}= p^{e(p,k)}\bigg(\prod_{{\rm prime} \ q:\ q\nmid
ac, \ q| L_k\atop p| (q-1), \ k+1\not\equiv 0\pmod {q}} \
q\bigg)\bigg(\prod_{{\rm prime} \ q:\  q\nmid a,\atop q| c, \ p|
(q-1)}q\bigg),
\end{align*}
 where
\begin{align*}
e(p,k):={\left\{
  \begin{array}{rl}
v_p(c), \quad&\text{if} \ v_{p}(k+1)\geq v_p(L_k),\\
v_p(cL_k), \quad&\text{if} \ v_p(k+1)< v_p(L_k).
 \end{array}
\right.}
\end{align*}
}

\begin{proof}
Similarly to the proof of Lemma 3.6, it is enough to
show that $p^{e(p,k)}$ is the smallest period of
$\sum_{e=1}^{v_p(cL_k)}f_e(n)$ by (3.3). We divide the proof into the
following two cases.

{\sc Case 1.} $v_p(k+1)\ge v_p(L_k)$. As in the proof of Subcase 2.2
 in Lemma 3.6, since $b'+a'n\equiv b'+a'(n+(k+1)c)\pmod{p^{v_p(cL_k)}}$, we can
obtain that $p^{v_p(c)}$ is a period of
$\sum_{e=1}^{v_p(cL_k)}f_e(n)$. If $v_p(c)=0$, it is complete. If
$v_p(c)\ge 1$, then choosing a positive integer $n_0$ such that
$v_p(b'+a'n_0)=v_p(c)$, we can show that $p^{v_p(c)-1}$ is not the
period of $\sum_{e=1}^{v_p(cL_k)}f_e(n)$ using a similar method as
in  the proof of Case 2 in Lemma 3.6.

{\sc Case 2.} $v_p(k+1)< v_p(L_k)$. Using  the same way as the proof
of Case 3 in  Lemma 3.6, one can easily check that
$p^{v_p(cL_k)-1}$ is not the period of
$\sum_{e=1}^{v_p(cL_k)}f_e(n)$.  The proof of Lemma 3.7 is complete.
\end{proof}

\section{\bf Proof of Theorem 1.3 and examples}

In this section, we first use the results presented in the previous
section to show Theorem 1.3. \\
\\
{\it Proof of Theorem 1.3.} By Theorem 1.2, we know that $g_{k,
\varphi}$ is periodic and $P_{k, \varphi}|cL_k$. To determine the
exact value of $P_{k, \varphi}$, it is sufficient to determine the
$p$-adic valuation of $P_{k, \varphi}$ for each prime $p$. By
Lemma 3.3, we have $P_{k, \varphi}={\rm lcm}_{{\rm prime}\ p\le
cL_k} \{P_{p, k,\varphi}\}$. So  it is enough to compute $\max_{{\rm
prime}\ q\le cL_k} \{v_p(P_{q, k, \varphi})\}$ for each prime $p$.
We consider the following four cases.

{\sc Case 1.} $p\nmid cL_k$.  Since $P_{k, \varphi}|cL_k$, it is
clear that $v_p(P_{k, \varphi})=v_p(cL_k)=0$.

{\sc Case 2.} $p| cL_k$ and $p|a'$. Observe from Lemmas 3.4-3.7 that
$v_p(P_{q, k, \varphi})=0$ for each prime $q\le cL_k$. So we have
$v_p(P_{k, \varphi})=0$.

{\sc Case 3.} $p=2$. From the discussion in Case 1 and Case 2, we
know that $v_2(P_{k,\varphi})=0$ if $2|cL_k$ and $2| a'$ or if
$2\nmid cL_k$. It remains to consider the case $2|cL_k$ and $2\nmid
a'$. By Lemmas 3.4-3.7, we know that $v_2(P_{p, k, \varphi})=0$ for
all odd primes $p$. So we only need to compute
$v_2(P_{2,k,\varphi})$.  We now distinguish the following four
subcases.

{\sc Subcase 3.1.} $2\nmid a$ and $v_2(cL_k)=1$. In this case, by
Lemma 3.6, one has $v_2(P_{k, \varphi})=0=v_2(cL_k)-1$.

{\sc Subcase 3.2.} $2\nmid a$, $v_2(cL_k)\ge 2$ and $v_2(k+1)\ge
v_2(L_k)=1$, or $2\nmid a', 2|d$ and $v_2(k+1)\ge v_2(L_k)=1$. Since
$v_2(k+1)\ge v_2(L_k)$ and $v_2(L_k)=1$, we get $k=3$. Thus by
Lemmas 3.6 and 3.7, we have that if $k=3$, $2\nmid a$ and
$v_2(c)=v_2(cL_k)-v_2(L_k)\ge 2-1=1$, or if $k=3$, $2\nmid a'$ and
$2|d$, then $v_2(P_{k, \varphi})=v_2(c)=v_2(cL_k)-1$.

{\sc Subcase 3.3.} $v_2(k+1)\ge v_2(L_k)\ge 2$.  Using Lemmas 3.6
and 3.7, we obtain that $v_2(P_{k,
\varphi})=v_2(c)=v_2(cL_k)-v_2(L_k)$.

{\sc Subcase 3.4.} $2\nmid a, v_2(L_k)=0$ and $v_2(cL_k)\ge 2$, or
$2\nmid a', 2|d$ and $v_2(L_k)=0$,  or $2\nmid a$, $v_2(cL_k)\ge 2$
and $v_2(k+1)<v_2(L_k)$, or $2\nmid a'$, $2|d$ and
$v_2(k+1)<v_2(L_k)$. If $2\nmid a, v_2(L_k)=0$ and $v_2(cL_k)\ge 2$,
or if  $2\nmid a', 2|d$ and $v_2(L_k)=0$, we get $v_2(P_{2,
k,\varphi})=v_2(c)=v_2(cL_k)$.  So  we have $v_2(P_{k,
\varphi})=v_2(cL_k)$ in this case.

Combining all the above information on $v_2(P_{k,\varphi})$, we have
\begin{align*}
v_2(P_{k, \varphi})={\left\{
\begin{array}{rl}
0, &\text{if} \ 2| a', \\
v_2(cL_k)-v_2(L_k), &\text{if} \ 2\nmid a' \ {\rm and}
\ v_2(k+1)\ge v_2(L_k)\ge 2, \\
v_2(cL_k)-1, &\text{if} \  2\nmid a \ \text{and}\ v_2(cL_k)=1,
\text{or} \ k=3, 2\nmid a\ {\rm and} \ 2|c,\\ &\quad  {\rm or} \
k=3, 2\nmid a'\ {\rm and}\ 2|d,\\
v_2(cL_k), &\text{otherwise}.
 \end{array}
\right.}(4.1)
\end{align*}

{\sc Case 4.} $p\ne 2, p|cL_k$ and $p\nmid a'$. Note that $2|(p-1)$
for each odd prime $p$. Evidently, if $2\nmid cL_k$, then $k=1$. So
there is no odd prime $p$ so that $p| L_k$ if $2\nmid cL_k$. Thus by
Lemmas 3.4-3.7, for all odd prime factors $p$ of $cL_k$, we obtain
that $ v_p(P_{2, k, \varphi})=1$ except that
 {either $p\nmid ac$}, $p|L_k$ and $k+1\equiv0\pmod p$ {or} $p|d$, in
which case $v_p(P_{2, k, \varphi})=0$. On the
other hand, for all odd primes $q$ such that $q\ne p$ and $q\le
cL_k$, we have by Lemmas 3.4-3.7 that $v_p(P_{q, k, \varphi})=0$
if $p|d$ or if $p\nmid ac$, $p|L_k$ and $k+1\equiv0\pmod p$, and $v_p(P_{q, k,
\varphi})\le 1$ otherwise. Hence $v_p(P_{2,k,\varphi})\ge v_p(P_{q,
k, \varphi})$ for all odd primes $q$ such that $q\ne p$ and $q\le
cL_k$.  Therefore we deduce immediately that
$$
v_p(P_{k, \varphi})=\max_{{\rm prime}\ q\le cL_k} \{v_p(P_{q, k,
\varphi})\}=\max( v_p(P_{2, k, \varphi}), v_p(P_{p, k,
\varphi})).\eqno(4.2)
$$
Using Lemmas 3.6 and 3.7 to compute $v_p(P_{p,k,\varphi})$, we get
that
\begin{align*}
v_p(P_{p, k, \varphi})={\left\{
  \begin{array}{rl}
0, &\text{if} \  v_p(cL_k)=1\ \mbox{and}\ p\nmid d,\\
v_p(c), &\text{if} \ v_{p}(k+1)\ge v_p(L_k), p\nmid d\
\mbox{and}\ v_p(cL_k)\ge 2,\\
&\ \mbox{or\ if}\ v_{p}(k+1)\ge v_p(L_k)\ \mbox{and}\ p|d,\\
v_p(cL_k), &\text{if} \ v_{p}(k+1)< v_p(L_k), p\nmid d\ \mbox{and}\
v_p(cL_k)\ge 2,\\
&\  \mbox{or\ if}\ \ v_{p}(k+1)< v_p(L_k)\ \mbox{and}\ p|d.
 \end{array}
\right.}\quad\quad\quad \quad\quad(4.3)
\end{align*}

For all the primes $p$ such that $p\nmid d$ and $v_p(cL_k)=1$, we have
by the above discussion  that $v_p(P_{2,
k, \varphi})=0$ only if $v_p(c)=0$, $v_p(L_k)=1$ and
$k+1\equiv0\pmod p$. Equivalently, $v_p(P_{2, k,
\varphi})=v_p(c)=v_p(cL_k)-v_p(L_k)$ if $v_p(k+1)\ge v_p(L_k)
\ge 1$, and $v_p(P_{2, k, \varphi})=v_p(cL_k)$ otherwise.
Therefore, for all the odd primes $p$ satisfying $p|cL_k$ and
$p\nmid a'$, we derive from (4.2) and (4.3) that
$$
v_p(P_{k, \varphi})={\left\{
  \begin{array}{rl}
v_p(cL_k)-v_p(L_k), &\text{if}\ v_{p}(k+1)\geq v_p(L_k)\ge 1,\\
v_p(cL_k), &\text{otherwise}.
 \end{array}
\right.}\eqno(4.4)
$$

Now putting all the above cases together, we get
\begin{align*}
P_{k, \varphi}&=2^{v_2(P_{k, \varphi})}\bigg(\prod_{{\rm prime}\ p:\
p\ne2 \atop p|a',\  p|cL_k}p^{v_p(P_{k,
\varphi})}\bigg)\bigg(\prod_{{\rm prime}\ p:\  p \ne 2 \atop p\nmid
a',\  p|cL_k}p^{v_p(P_{k,\varphi})}\bigg)\\
&=\frac{cL_k}{2^{v_2(cL_k)-v_2(P_{k, \varphi})}\bigg(\prod_{{\rm
prime}\ p:\ p\ne 2\atop  p|a',\ p|cL_k}p^{v_p(cL_k)-v_p(P_{k,
\varphi})}\bigg)\bigg(\prod_{{\rm prime}\ p:\ p\ne 2
\atop p\nmid a',\ p|cL_k}p^{v_p(cL_k)-v_p(P_{k,\varphi})}\bigg)}\\
&=\frac{cL_k}{2^{v_2(cL_k)-v_2(P_{k,\varphi})}\bigg(\prod_{{\rm
prime}\ q:\ q\ne 2,\ q|a'}q^{v_q(cL_k)}\bigg)\bigg(\prod_{{\rm
prime}\ q:\ q\ne 2\atop  q\nmid
a',\ q|cL_k}q^{v_q(cL_k)-v_q(P_{k,\varphi})}\bigg)}\\
&=\frac{cL_k}{2^{\delta_{2, k,\varphi}}\bigg(\prod_{{\rm prime}\
q|a'}q^{v_q(cL_k)}\bigg)\bigg(\prod_{{\rm prime}\ q:\ q\ne 2 \atop
q\nmid a',\  q|cL_k}q^{v_q(cL_k)-v_q(P_{k,\varphi})}\bigg)}, \ \ \ \
\ \ \ \ \ \  \ \ \ \ \  \ \ \ \ \ (4.5)
\end{align*}
where
\begin{align*}
\delta_{2, k,\varphi}:={\left\{
\begin{array}{rl}
v_2(cL_k)-v_2(P_{k, \varphi }), &\text{\rm if} \ 2\nmid a',\\
0, &\text{\rm if} \ 2|a'.
 \end{array}
\right.}
\end{align*}
It then follows from (4.1) that
\begin{align*}
\delta_{2, k,\varphi}={\left\{
\begin{array}{rl}
v_2(L_k), &\text{\rm if} \ 2\nmid a' \ {\rm and} \ v_2(k+1)\ge v_2(L_k)\ge 2, \\
1, &\text{\rm if} \  2\nmid a \ {\rm and}\ v_2(cL_k)=1, {\rm or} \
k=3,
2\nmid a\ {\rm and} \ 2|c,  {\rm or} \ k=3, 2\nmid a'\ {\rm and}\ 2|d,\\
0, &\text{\rm otherwise,}
 \end{array}
\right.}
\end{align*}
which implies that  $\eta_{2,k, a', c}=2^{\delta_{2,k,\varphi}}$.
Hence by (1.2) we get
$$Q_{k, a', c}=\frac{cL_k}{2^{\delta_{2, k,\varphi}}\prod_{{\rm
prime}\ q|a'}q^{v_q(cL_k)}}. \eqno(4.6)$$

Since there is at most one odd prime $p\le k$ such that $v_p(k+1)\ge
v_p(L_k)\ge 1$ (see \cite{[FK]}), we derive from (4.4) that
\begin{align*}
\prod_{{\rm prime}\ q:\ q \ne 2,\atop  q\nmid a',\
q|cL_k}q^{v_q(cL_k)-v_q(P_{k,\varphi})}={\left\{
  \begin{array}{rl}
p^{v_p(L_k)}, &\text{if} \  v_p(k+1)\geq v_p(L_k)\ge 1
\ \mbox{for an odd prime}\ p\nmid a',\\
1, &\text{otherwise.}
 \end{array}
\right.}
\end{align*}
Thus it follows from (4.5) and (4.6) that $P_{k,\varphi}$ is equal
to $Q_{k, a', c}$ except that $v_p(k+1)\ge v_p(L_k)\ge 1$
for at most one odd prime $p\nmid a'$, in which case
$P_{k, \varphi}$ equals  $\frac{Q_{k, a', c}}{p^{v_p(L_k)}}$.

The proof of Theorem 1.3 is complete.
\hfill$\square$\\

Now we give some examples to illustrate Theorem 1.3.\\
\\
{\bf Example 4.1.} Let $ a\ge 1, b\ge 0$ and $c\ge 1$ be integers,
and let $a':=a/\gcd(a, b)$ be odd. Let $k=2^t-1$, where $t\in
\mathbb{N}$ and  $t\ge 3$.  Since $v_2(k+1)=t>v_2(L_k)=t-1\ge 2$, we
obtain by Theorem 1.3 that $\eta_{2,k,a',c}=2^{v_2(L_k)}$. On the
other hand, there is no odd prime $p$ satisfying $v_p(k+1)\ge v_p(L_k)\ge 1$.
Thus we have $$ P_{k, \varphi}=\frac{cL_k}{2^{v_2(L_k)}
\prod_{{\rm prime} \ q|a'}q^{v_q(cL_k)}}.
$$\\
\\
{\bf Example 4.2.}    Let $ a\ge 1, b\ge 0$ and $c\ge 1$
be integers, and let  $a':=a/{\rm gcd}(a, b)$.  Let $p$ be any given
odd prime with $p\nmid a'$, and let $k=p^{\alpha}-1$ for some
integer $\alpha\ge 2$. Since $k=p^{\alpha}-1>3$ and
$v_2(k+1)=v_2(p^{\alpha})=0$, we have $\eta_{2,k, a',c}=1$. The odd
prime $p$ satisfies that $v_p(k+1)=\alpha>\alpha-1=v_p(L_k)\ge 1$.
Hence we get by Theorem 1.3 that
$$
P_{k, \varphi}=\frac{cL_k}{p^{v_p(L_k)}\prod_{{\rm prime} \
q|a'}q^{v_q(cL_k)}}.
$$\\
\\
{\bf Example 4.3.}   Let $ a\ge 1, b\ge 0$ and $c\ge 1$ be integers,
and let  $a':=a/{\rm gcd}(a, b)$.  If $k$ is an integer of the form
$35^{\alpha}-1$ with  $\alpha\ge 2$ and $\alpha\in \mathbb{N}$, then
$$
P_{k, \varphi}=\frac{cL_k}{\prod_{{\rm prime} \ q|a'}q^{v_q(cL_k)}}.
\eqno(4.7)
$$
Actually, since $35^{\alpha}-1>3$ and $v_2(35^{\alpha})=0<v_2(L_k)$,
we obtain by Theorem 1.3 that $\eta_{2, k,a',c}=1$. On the other
hand, we have that  $v_5(k+1)=\alpha<v_5(L_k)$,
$v_7(k+1)=\alpha<v_7(L_k)$ and $v_q(k+1)=0$ for any other odd prime
$q$. Hence we get $P_{k,\varphi}$ as in (4.7).

Furthermore, if $a|b$ or $a'$ is a prime
greater than $cL_k$,  then there is no prime factor of $a'$ dividing
$cL_k$. Therefore  for any $k=35^{\alpha}-1$ with $\alpha\ge 2$ and
$\alpha\in \mathbb{N}$, one has $P_{k, \varphi}=cL_k$.\\

Finally, by Theorem 1.3, we only need  to compute the first $P_{k,
\varphi}$ values of $g_{k,\varphi}$ so that we can estimate the
difference between $\prod_{0\le i\le k} \varphi(b+a(n+ic))$ and
$\varphi({\rm lcm}_{0\le i\le k}\{b+a(n+ic)\})$ for large $n$. In
other words, we have
$$\min_{1\le m\le P_{k,\varphi}}\{g_{k,\varphi}(m)\}
\le \frac{\prod_{0\le i\le k} \varphi(b+a(n+ic))} {\varphi({\rm
lcm}_{0\le i\le k}\{b+a(n+ic)\})}=g_{k,\varphi}\big(\langle
n\rangle_{P_{k,\varphi}}\big)\le \max_{1\le m\le
P_{k,\varphi}}\{g_{k,\varphi}(m)\},$$ where $\langle
n\rangle_{P_{k,\varphi}}$ means the integer between $1$ and
$P_{k,\varphi}$ such that $n\equiv\langle
n\rangle_{P_{k,\varphi}}\pmod{P_{k,\varphi}}$.

On the other hand, estimating the difference between $\prod_{0\le
i\le k}\varphi(b+a(n+ic))$ and ${\rm lcm}_{0\le i\le
k}\{\varphi(b+a(n+ic))\}$ is also an interesting problem. For this
purpose, we define the arithmetic function $G_{k, \varphi}$ for any
positive integer $n$ by
$$
G_{k, \varphi}(n):=\frac{\prod_{i=0}^k \varphi(b+a(n+ic))} {{\rm lcm}_{0\le i\le k}
\{\varphi(b+a(n+ic))\}}.
$$
Unfortunately, $G_{k, \varphi}$ may not be periodic. For instance,
taking $a=1, b=0$ and $c=1$, then the arithmetic function $\bar
G_{k, \varphi}$ defined by $\bar G_{k,
\varphi}(n):=\frac{\prod_{i=0}^n \varphi(n+i)}{{\rm lcm}_{0\le i\le
k}\{\varphi(n+i)\}}$ for $n\in \mathbb{N}^*$ is not periodic.
Indeed, for any given positive integer $M$, we can always choose a
prime $p>M$ since there are infinitely many primes. By Dirichlet's
theorem, we know that there exists a positive integer $m$ such that
the term $mp^2+1$ is a prime in the arithmetic progression
$\{np^2+1\}_{n\in \mathbb{N}^*}$. Letting $n_0=mp^2$ gives us that
$p|\varphi(n_0)$ and $p| \varphi(n_0+1)=\varphi(mp^2+1)=mp^2$. Thus
$p|\bar G_{k, \varphi}(n_0)$ and $\bar G_{k, \varphi}(n_0)\ge p>M$.
That is, $\bar G_{k, \varphi}$ is unbounded, which implies that
$\bar G_{k, \varphi}$ is not periodic. Applying Theorem 1.3, we can
give a nontrivial upper bound about the integer ${\rm lcm}_{0\le
i\le k}\{\varphi(b+a(n+ic))\}$ as follows.\\
\\
\noindent{\bf Proposition 4.4.} {\it Let $k\ge 1, a\ge 1, b\ge 0$
and $c\ge 1$ be integers. Then for any positive integer $n$, we have
\begin{align*}
{\rm lcm}_{0\le i\le k}\{\varphi(b+a(n+ic))\} \le
\frac{\prod_{i=0}^k\varphi(b+a(n+ic))}{g_{k,\varphi}\big(\langle
n\rangle_{P_{k,\varphi}}\big)}
\end{align*}
with $\langle n\rangle_{P_{k,\varphi}}$ being defined as above.}
\begin{proof}
For each $0\le i\le k$, since $\varphi$ is multiplicative and $b+a(n+ic)|
{\rm lcm}_{0\le j\le k}\{b+a(n+jc)\}$, we have
$$
\varphi(b+a(n+ic))\big|
\varphi({\rm lcm}_{0\le j\le k}\{b+a(n+jc)\}).
$$
So we get ${\rm
lcm}_{0\le i\le k}\{\varphi(b+a(n+ic))\}\big| \varphi({\rm
lcm}_{0\le i\le k}\{b+a(n+ic)\})$. Thereby
$$g_{k,\varphi}\big(\langle n\rangle_{P_{k, \varphi}}\big)=g_{k,\varphi}(n)\le
\frac{\prod_{0\le i\le k}\varphi(b+a(n+ic))}{{\rm lcm}_{0\le i\le
k}\{\varphi(b+a(n+ic))\}}
$$
for any positive integer $n$. The desired result then follows immediately.
\end{proof}

Theorem 1.3 answers the second part of Problem 1.1 for the Euler phi
function. However, the smallest period problem is still kept open for all other
multiplicative functions $f$ with $f(n)\ne 0$ for all positive
integers $n$. For example, if one picks $f=\sigma _{\alpha}$ with
$\sigma _{\alpha}(n):=\sum_{d|n\atop d\ge 1}d^{\alpha }$ for $\alpha
\in \mathbb {N}$, then what is the smallest period of $g_{k, f}$?
If $f=\xi _{\varepsilon}$ with
$\xi _{\varepsilon}(n):=n^\varepsilon$ for $\varepsilon \in \mathbb{R}$,
then what is the smallest period of $g_{k, f}$?\\ 

\begin{center}
{\bf Acknowledgements}
\end{center}
The authors are grateful to the anonymous referees for careful
reading of the manuscript and for helpful comments and suggestions.

\end{document}